\documentclass[reqno, letterpaper, final, sort&compress]{amsart}

\usepackage[margin=3cm]{geometry}
\usepackage[utf8]{inputenc}  
\usepackage[T1]{fontenc}

\usepackage{relsize}  
\usepackage{amssymb}  
\usepackage{mathtools}
\usepackage{amsfonts,amscd,amsthm}
\usepackage{graphicx} 
\graphicspath{{pics/}}

\usepackage[dvipsnames]{xcolor}  
\usepackage[section]{placeins}

\usepackage{fontawesome5}
\newcommand{\orcid}[1]{\href{https://orcid.org/#1}{\textcolor[HTML]{A6CE39}{\faOrcid}}}

\usepackage{tikz}
\usetikzlibrary{positioning}
\usetikzlibrary{backgrounds}
\usetikzlibrary{calc}

\definecolor{zzttqq}{rgb}{0.6,0.2,0}
\definecolor{mBlue}{rgb}{0.3,0.3,1}

\usepackage{setspace}
\usepackage{booktabs}
\newcommand{\otoprule}{\midrule[\heavyrulewidth]}%
\usepackage{hyperref}
\hypersetup{plainpages=false, linktocpage=true, colorlinks=false, pdfborderstyle={/S/U/W 1},
            breaklinks=true, linkcolor=black, menucolor=black,
			urlcolor=black, citecolor=black}

\usepackage{upquote}
\usepackage{fancyvrb}

\RecustomVerbatimCommand{\VerbatimInput}{VerbatimInput}%
{fontsize=\footnotesize,
 frame=lines,  
 framesep=1em, 
 rulecolor=\color{Gray},
 label=\fbox{\color{Black}\$HOME/Documents/conexp/runConexp.sh},
 labelposition=topline,
 commentchar=*        
}

\usepackage{enumitem}

\newtheorem{lemma}{Lemma}

\newtheorem*{corollary*}{Corollary}
\newtheorem*{lemma*}{Lemma}

\theoremstyle{definition}

\newtheorem{example}[lemma]{Example}

\providecommand{\ConExp}{ConExp}
\providecommand{\Java}{\textsc{Java}}
\providecommand{\Oracle}{\textsc{Oracle}}

\newcommand{\K}{\mathbb{K}}

\definecolor{mBlue}{HTML}{006AD7}%

\usepackage{mleftright}

\newcommand{\set}[1]{\ensuremath{\mleft\{ #1 \mright\}}}

\newcommand{\apply}[1]{\ensuremath{\mleft( #1 \mright)}}

\newcommand{\defeq}{\coloneqq}
\newcommand{\eqdef}{\mathrel{\mathopen={\mathclose:}}}

\providecommand{\GoodNews}{\colorbox{green!20!white}{\color{black}Good News:}}
\providecommand{\green}[1][\checkmark]{\colorbox{green!20!white}{\color{black}#1}}
\providecommand{\yellow}[1][\checkmark]{\colorbox{yellow!20!white}{\color{black}#1}}

\title{FCA using the Concept Explorer in 2024}
\thanks{\textit{E-mail addresses:} \texttt{edith.vargas@itam.mx}, \texttt{andreas.wachtel@itam.mx}}
\thanks{Departamento Académico de Matemáticas, ITAM, México}
\thanks{This work was supported  by the Asociación Mexicana de Cultura A.C}

\author[Edith Vargas-Garc\'\i{}a]{\orcid{0000-0001-9677-9087}~Edith Vargas-Garc\'\i{}a \ }

\author[Andreas Wachtel]{\ Andreas Wachtel~\orcid{0000-0002-6317-787X}}

\newcommand{\codebckg}[1]{{\setlength{\fboxsep}{1pt}\colorbox{gray!20!white}{\color{black}\tt#1}}}

\begin{document}

\begin{abstract}
	In this note we give a very short introduction to Formal Concept Analysis,
accompanied by an example in order to build concept lattices from a context.
We build the lattice using the \Java-based software \emph{Concept Explorer} (\ConExp) in a recent version of Linux.

Installing an appropriate \Java\ version is necessary,
because \ConExp\ was developed some time ago using 
a \textsc{Sun} \Java\ version,
which is not open-source.
As a result,  it has been observed that \ConExp\ will not build a lattice when started with an
  open-source  \Java\ version.
Therefore,  we also sketch the procedure 
we followed to install an appropriate \Java\ version
which makes \ConExp\ work again, \emph{i.e.}, to ``build lattices again''.
We also show how to start \ConExp\ with a 32\,bit \Java\ version,
which requires a few additional libraries.

\end{abstract}

\maketitle

2010 Mathematics Subject Classification: 06D50, 06A06, 06B05

\textbf{Keywords.}  Lattices, Partial Order, Concept Explorer, FCA

\smallskip
\subsection*{Declaration}
The purpose of this note is to 
make Formal Concept Analysis a bit more visible and accessible, give some references to theory
and to show that the \Java\ software called the \emph{Concept Explorer} (\ConExp)
can still be used to build lattices.
We were asked to write this note, because we managed to make \ConExp\ build lattices again and thus make it usable again in 2024,
almost 18 years after its last publication date in 2006.

This note is organized as follows.
In Section~\ref{secPreliminaries} we introduce basic notions of Formal Concept Analysis.
Section~\ref{secConExp} 
illustrates how to use of the Concept Explorer on a small example
and
gives references to its official documentation.
Finally, Sections~\ref{secInstallation} and \ref{secInstallation32} list compatible \Java\ versions, available today, 
and describe the procedure we used to install a working Concept Explorer.

\bigskip
\section{Basic Notions}\label{secPreliminaries}

Formal concept analysis (FCA) is a theoretical framework for data analysis based on the notion of a \emph{formal concept},
it harnesses the powers of general abstract Galois theory and the structure theory of complete lattices. 
FCA was born around 1980, when a research group in Darmstadt, Germany, headed by Rudolf Wille, began to develop a framework for lattice theory applications, \cite{ganter1997applied}:
For sets~$G$ and~$M$ and any binary relation $I\subseteq G\times M$ between~$G$ and~$M$ the triple $\mathbb{K}=\apply{G,M,I}$ is called a \textbf{formal context}, and~$I$ its \textbf{incidence relation}. 
The elements of~$G$ are called \textbf{objects} and those of~$M$ are \textbf{attributes}. 
When~$G$ and~$M$ are finite, a  formal context can be represented by a cross table. 
The elements on the left-hand side (row labels) are the objects, the elements at the top (column labels) are the attributes, and the incidence relation is represented by the crosses.

\clearpage
\begin{example}\label{Exa1}
\textbf{Planets of our solar system.} 
The following context with $9$ objects and $7$ attributes, shows some characteristics of the planets.

\begin{table}[h!]
\centering
\providecommand{\TX}{\ensuremath{\times}}
\providecommand{\footns}[1]{{\footnotesize #1}}

\resizebox{0.95\textwidth}{!}{%
\begin{tabular}{|l||c|c|c|c|c|c|c|c|c|}\hline
& \footns{Size small} & \footns{Size medium} & \footns{Size large} & \footns{Near to sun} & \footns{Far from sun} & \footns{Moon yes} & \footns{Moon no}\\
& \footns{(small)} & \footns{(medium)} & \footns{(large)} & \footns{(near sun)} & \footns{(far from sun)} & \footns{(moon)} & \footns{(no moon)}
	\\ \hline\hline
\footns{Mercury (Me)} & \TX &     &     & \TX &     &     & \TX \\ \hline
\footns{Venus (V)}    & \TX &     &     & \TX &     &     & \TX \\ \hline
\footns{Earth (E)}    & \TX &     &     & \TX &     & \TX &     \\ \hline
\footns{Mars (Ma)}    & \TX &     &     & \TX &     & \TX &     \\ \hline
\footns{Jupiter (J)}  &     &     & \TX &     & \TX & \TX &     \\ \hline
\footns{Saturn (S)}   &     &     & \TX &     & \TX & \TX &     \\ \hline
\footns{Uranus (U)}   &     & \TX &     &     & \TX & \TX &     \\ \hline
\footns{Neptune (N)}  &     & \TX &     &     & \TX & \TX &     \\ \hline
\footns{Pluto (P)}    & \TX &     &     &     & \TX & \TX &     \\ \hline
    \end{tabular}
}

\caption{\rule{0pt}{3ex}Planets and their attributes of our solar system, taken from \cite[p.65]{davey}.}\label{tabContextPsychometric}
\end{table}
\end{example}

In order to define the formal concepts we will use the following two derivation operators, which are a central notion in FCA.
We write that an element $g\in G$ is related to an attribute $m\in M$ using the relation $I$ as $gIm$ and say: the object $g$ has the attribute $m$.  
Let $A\subseteq G$ and $B\subseteq M,$ then
the set of attributes common to all the objects in $A$ is
\begin{align*}
A'&\defeq \set{m\in M\mid \forall g\in A\colon gIm },
\end{align*}
and the set of objects possessing the attributes in $B$ is
\begin{align*}
B'&\defeq \set{g\in G \mid \forall m\in B\colon gIm }.
\end{align*}
The pair $({'},{'})$ of derivation operators establishes a \textbf{Galois connection} between the power set lattices on $G$ and $M,$ induced by $I,$ see also 
\cite[Chapter~V, Theorem~2.3]{Hungerford2011}. 

A \textbf{formal concept} is a pair $(A,B)$ where the \textbf{extent}
$A\subseteq G$ and the \textbf{intent} $B\subseteq M$ are sets 
with 
$A=B'$ and $B=A',$ \cite[page 18]{ganter}.
Intents of the form
$B=\set{g}'$ with $g\in G$ are called \textbf{object intents} and are
written as $g',$ for short; dually extents $A=\set{m}'\eqdef m'$ with
$m\in M$ are referred to as \textbf{attribute extents}.

Given two concepts $(A_1,B_1)$ and $(A_2,B_2),$
using one of the derivations operators we have that $A_1\subseteq A_2 \implies A_2' \subseteq A_1',$ 
wherefore $B_2=A_2' \subseteq A_1'=B_1.$
Thus, we define an order between two concepts as follows.
We write 
$(A_1,B_1)\leq (A_2,B_2)$ and say 
that $(A_1,B_1)$ is a \textbf{subconcept} of $(A_2, B_2),$  if $A_1\subseteq A_2$ or equivalently $B_1\supseteq B_2.$
The set of all concepts of $\mathbb{K} =\apply{G, M, I}$ ordered in this way is denoted by $\underline{\mathfrak{B}}\apply{\mathbb{K}}$ and  called the \textbf{concept lattice} of $\K.$

The fundamental theorem of formal concept
analysis~\cite[Theorem~3]{ganter} states that  every complete
lattice is a concept lattice, up to isomorphism. 
In particular, 
if~$\mathfrak{L}=\apply{L, \leq}$ is a complete lattice, then
$\mathfrak{L}\cong\underline{\mathfrak{B}}(L,L,{\leq}).$

\begin{figure}[h!]%
\centering
\tikzstyle{object}=[align=center, text width=10pt, fill=white, very thin, draw=gray!50, inner sep=1.25pt, outer sep=1pt]
\tikzstyle{attribute}=[fill=white, very thin, draw=gray!50, inner sep=1.5pt, outer sep=2pt]
\tikzstyle{normal}=[circle, draw=black, fill=black,inner sep=2.75pt, outer sep=0.125pt]
\tikzstyle{smaller}=[circle, draw,very thin,fill=white,inner sep=1.75pt, outer sep=0.125pt]
\tikzset{ntJ/.style={path picture={\fill[white](path picture bounding box.east) rectangle (path picture bounding box.north west);}}}
\tikzset{ntJb/.style={path picture={\fill[mBlue](path picture bounding box.east) rectangle (path picture bounding box.north west);}}}
\tikzset{ntM/.style={path picture={\fill[white](path picture bounding box.east) rectangle (path picture bounding box.south west);}}}

\begin{singlespace}
\begin{tikzpicture}[scale=0.75, x=-1cm]
\node[smaller] (12) at (1, 10) {};
\node[smaller,ntJb,label={[attribute]90:{\tiny{moon}}}] (11) at (3, 8) {};
\node[smaller,ntJb,label={[attribute]90:{\tiny{small}}}] (10) at (-1, 8) {};
\node[smaller] (8) at (1, 6) {};
\node[smaller,ntJb, label={[attribute]90:{\tiny{far from sun}}}] (9) at (5, 6) {};
\node[smaller,ntJb, label={[attribute]90:{\tiny{near sun}}}] (7) at (-3, 6) {};
\node[normal,ntJ,label={[object]270:{\tiny{E Ma}}}] (4) at (-1, 4) {};
\node[normal,ntJ, label={[object]270:{\tiny{P}}}] (5) at (3, 4) {};
\node[smaller] (6) at (1, 2) {};
\node[normal,ntJb, label={[object]270:{\tiny Me V}}, label={[attribute]90:{\tiny{no moon}}}] (3) at (-5, 4) {};
\node[normal,ntJb, label={[object]270:{\tiny U N}}, label={[attribute]90:{\tiny{medium}}}] (2) at (7, 4) {};
\node[normal,ntJb, label={[object]270:{\tiny{J S}}}, label={[attribute]90:{\tiny{large}}}] (1) at (5, 4) {};

\scoped[on background layer]{%
\begin{scope}
\draw (12)--(11)--(8)--(10)--(12);
\draw (10)--(7);
\draw (11)--(9);
\draw (9)--(5)--(8);
\draw (7)--(4)--(8);
\draw (6)--(3);
\draw (6)--(2);
\draw (6)--(1);
\draw (6)--(4);
\draw (6)--(5);
\draw (9)--(2);
\draw (9)--(1);
\draw (7)--(3);
\end{scope}
}

\end{tikzpicture}
\end{singlespace}

\caption{The concept lattice of the planetary Context~\ref{tabContextPsychometric}, see also 
\cite[p.74]{davey}.}\label{fig1}
\end{figure}
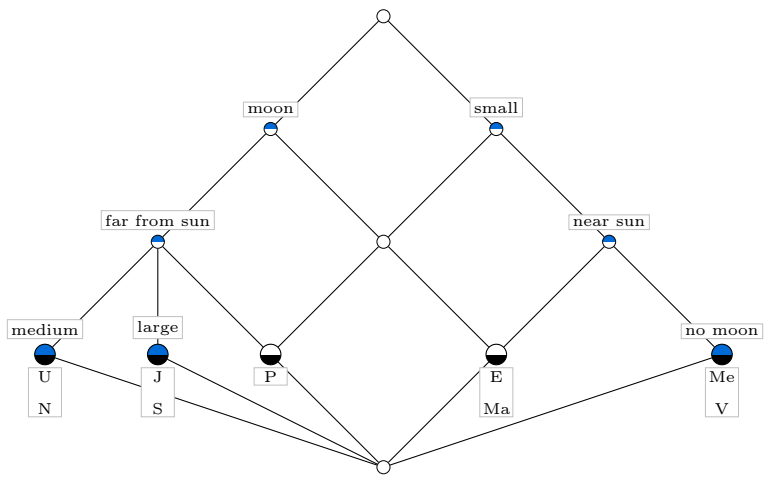

The concept lattice of Example \ref{Exa1} is depicted in Figure \ref{fig1}.
Each node represents a formal concept.
Object labels are depicted slightly below  and attribute labels slightly above nodes.
If a node has a blue upper semi-circle, then there is an attribute attached to this concept.
If a node has  a black lower semi-circle, then there is an object attached to this concept.
\clearpage

To give an example of a formal concept in Figure~\ref{fig1}
we pick
the central empty circle  which is associated to 
the formal concept $\apply{\set{\text{Earth, Mars, Pluto}}, \set{\text{size small, moon-yes}}}.$
The object- and attribute- sets of this concept 
can be read from the lattice as follows:
The objects in the extend of the concept are found ``walking down'' from the concept node to the bottom of the lattice, while attributes in the intent of the concept are found by ``walking down'' from the top of the lattice to   the concept node.
In the lattice ``walking down'' means that you are not allowed to return or go up and down.

\bigskip
\section{The Concept Explorer}\label{secConExp}

The Concept Explorer (\ConExp) implements many ideas needed for the study of data and research using Formal Concept Analysis.
\ConExp\ can (among other things) represent the structure of a finite formal context in form of a cross table
and can generate  a concept lattice.
This is not only a graph with nodes and edges connecting them, but an ordered structure (complete lattice) with a bottom and a top element. 
Each node of a lattice corresponds to a formal concept $\apply{O,A},$ where $O$ is a subset of the object set, $A$ is a subset of the attribute set.

Next, we give the basic steps to build a lattice of a given context in \ConExp.
More advanced steps for data analysis with \ConExp\ can be found in the complete documentation
 available at:
\begin{center}
	\url{https://conexp.sourceforge.net/users/documentation}.
\end{center}

\bigskip
If it is desired to use  \ConExp\ and to reproduce the figures while reading section \ref{secConExp},
then we recommend completing the installation steps in Section~\ref{secInstallation} before continuing to read this section.
\clearpage
\subsection{Basic Steps to use the Concept Explorer}\label{secConExpContext}
Starting \ConExp\ will be explained in Section~\ref{secStartingConExp}.
For the moment, suppose we just started \ConExp.
Then, the following window will appear.
The number of objects and attributes (initially 15) can be adjusted by changing the parameters in the left pane (where the mouse-pointer is).
\begin{center}
\includegraphics[scale=0.32]{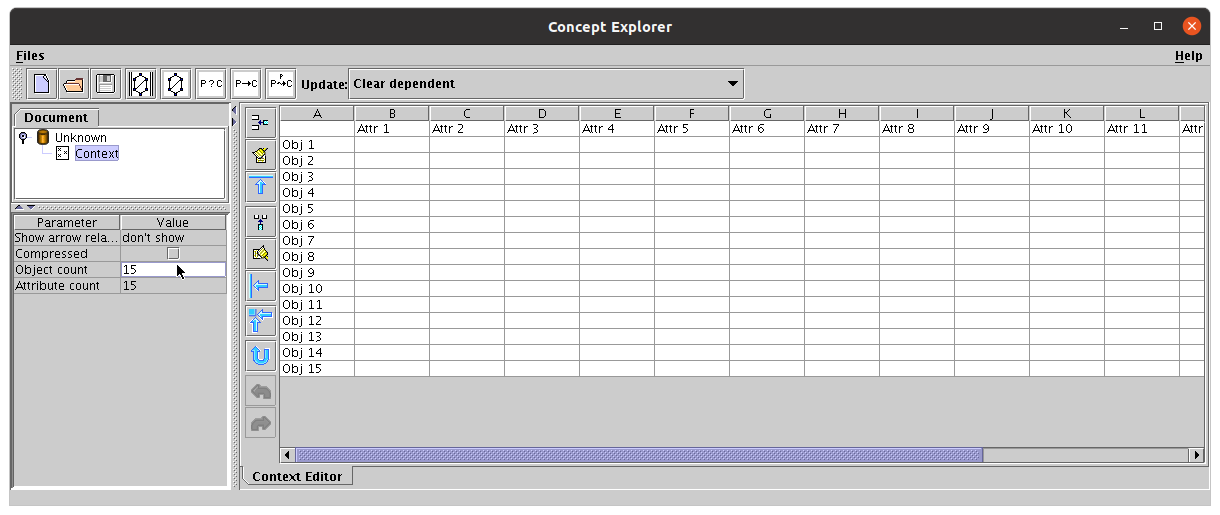}
\end{center}
		In the top left of the Concept Explorer window there are a few buttons.
	The first 3 stand for ``new'', ``open'' and ``save'' context.
The Context~\ref{tabContextPsychometric}
can be put into \ConExp,
by loading it from a context file \texttt{planets.cxt},
or by putting it in by hand.
To do that, the  object and attribute names can be changed
and the crosses can be inserted by hitting the ``x-key''
and to skip a cell one can use the ``dot-key''.
This is also  well explained in the documentation.
Once finished, after reducing the number of objects and attributes,
and after saving it to the context file \texttt{planets.cxt},
the Context~\ref{tabContextPsychometric} may look as follows: 
\begin{center}
	\includegraphics[scale=0.32]{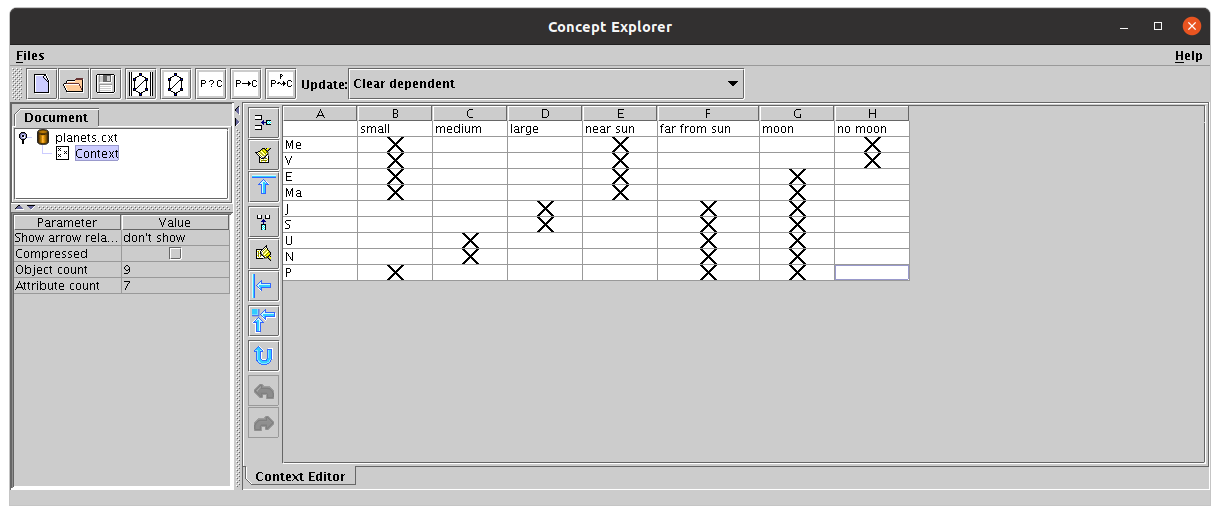}
\end{center}

Moreover, there are 2 buttons with a small lattice symbol 
\raisebox{-0.15\height}{\includegraphics[scale=0.5]{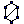}},
the left-one calculates the number of concepts, also the number of nodes in the concept lattice, and shows it in a small message window:
\begin{center}
	\includegraphics[scale=0.32]{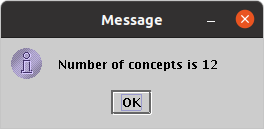}
\end{center}

\clearpage
\subsection{Displaying the concept lattice}\label{secConExpBuildLattice}
After putting in the context, as shown above,
clicking	the second  button with a small lattice symbol
\raisebox{-0.15\height}{\includegraphics[scale=0.5]{lattice}},
builds a concept lattice, provided  \ConExp\ is started properly (with a compatible \Java\ version, see Section~\ref{secInstallation}).
The result for Context~\ref{tabContextPsychometric} looks initially as follows,
the lattice is built,  but no labels are shown:

\begin{center}
	\includegraphics[scale=0.32]{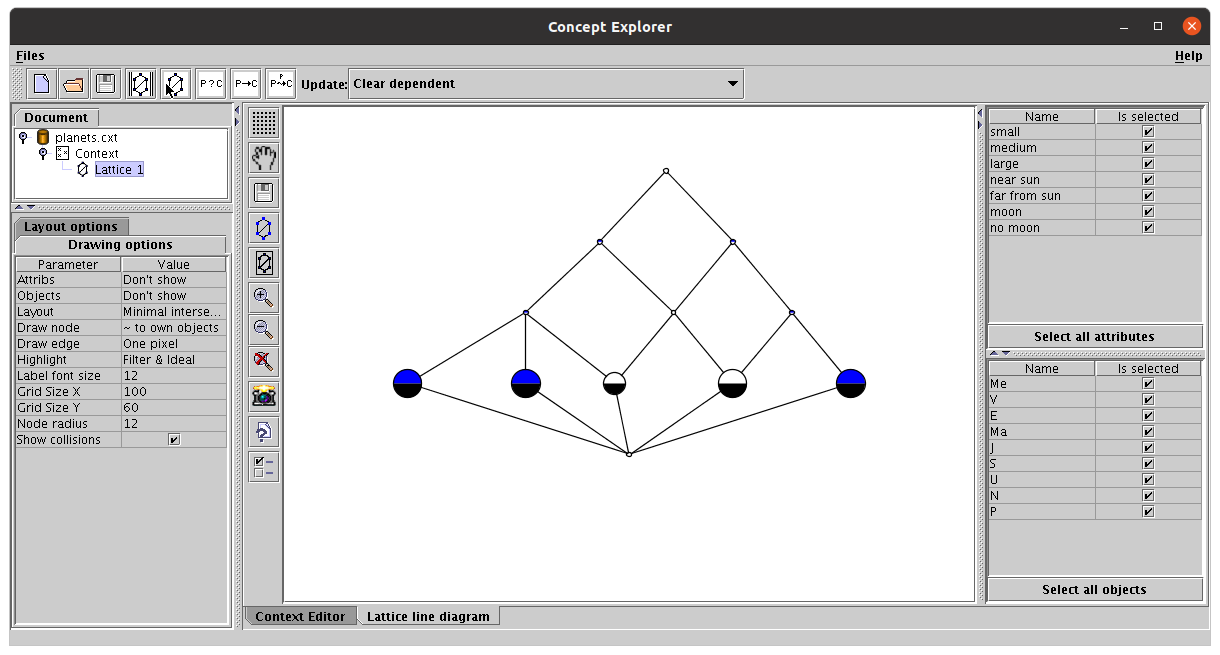}
\end{center}

After showing the labels of objects and attributes (see left pane under \emph{Drawing options}),
 hiding the right pane (using the small triangular buttons) and  adjusting the zoom,
 the result might look as follows.
 An awesome feature of \ConExp\ is that the nodes and labels can be moved by dragging them with the mouse.

\begin{center}
\includegraphics[scale=0.32]{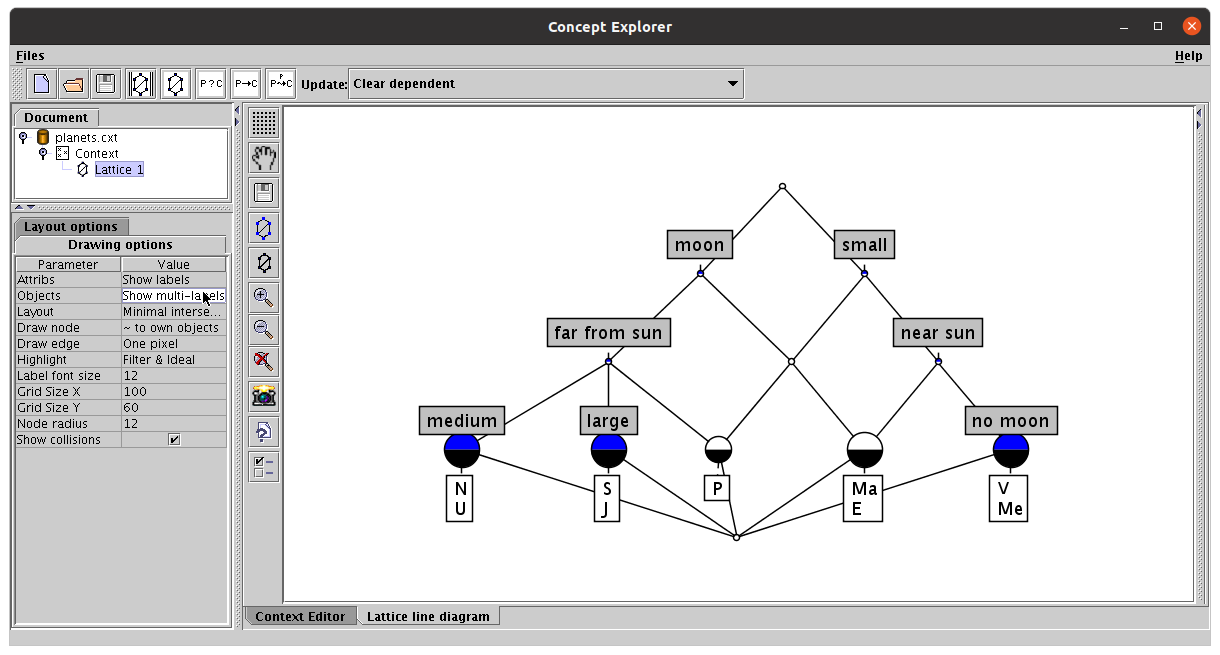}
\end{center}

\bigskip
\bigskip
It is beyond the scope of this note to discuss all FCA techniques.
However, an important technique that must be
mentioned is the use of ``implications'' between attributes in a formal context.

\clearpage
\subsection{Attribute implications}\label{secConExpImplications}
The implications are statements of the following
kind: \emph{Every object with the attributes $a, b, c, \ldots$ also has the attributes}
$x, y, z,\ldots$. 
Formally, an \textbf{implication between attributes}, is a
pair of subsets of the attribute set $M,$ denoted by $A \implies B.$
The set $A$ is the \textbf{premise} of the implication $A \implies B $, and $B$ is its \textbf{conclusion}. 
We say that $A \implies B$ holds in a context $\apply{G, M, I},$ if every object intent respects $A \implies B,$ that is, if each object that has all the attributes from $A$ also has all the
attributes from $B.$
Implications have been studied by Ganter \& Wille (since 1986),
more information  can be found in \cite[pp.\ 79]{ganter}.

In implications between attributes only the number of objects that satisfy the implication is important,
because it says  how many objects support or witness the implication.
In \ConExp\ clicking the button 
{%
\setlength{\fboxsep}{0pt}%
\setlength{\fboxrule}{0.25pt}%
\raisebox{-0.2\height}{\fbox{\includegraphics[scale=0.5]{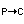}}}%
}
will calculate the implications and show them  in the following format:
\[
\text{ID } < \text{ Number of objects } > \text{ Premise } \implies \text{ Conclusion}.
\]
Showing the implications of Example~\ref{Exa1}  with \ConExp\ looks as follows:
\begin{center}
\includegraphics[scale=0.32]{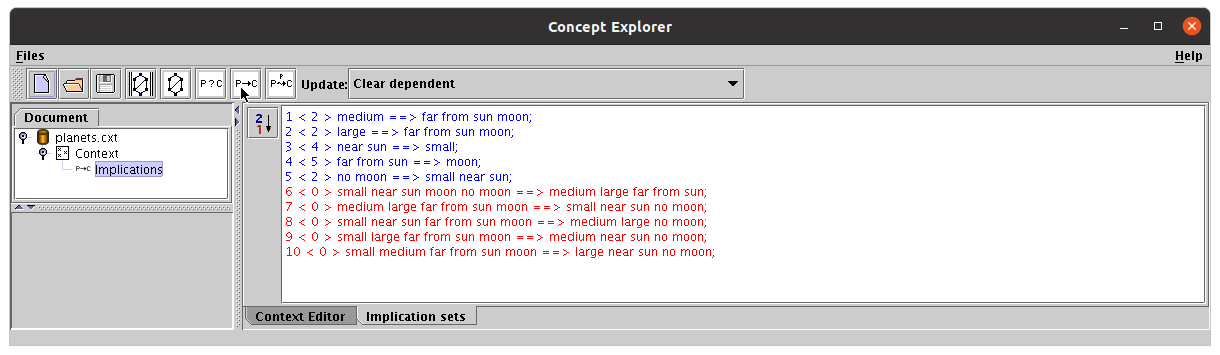}
\label{Implications}
\end{center}

Implications are displayed by \ConExp\ in either  blue or red,
where  blue indicates that there are objects in context which support the rule
and  red  indicates that there are no objects supporting the highlighted  rule. 

In the picture above,
the  implication with $\text{ID} = 1$ is:
\[
\text{1 } < \text{ 2 } > \text{ medium } \implies \text{ far from sun, moon}.
\]
It states that \textbf{all} objects (in this case 2 planets)  that have the attribute (medium),
also have the attribute (moon) and the attribute (far from sun).
If we look at Context~\ref{tabContextPsychometric} or the lattice in Figure~\ref{fig1}, then these 
objects are the planets Uranus and Neptune.
Furthermore, both the lattice and the context show  that \textbf{all} medium sized planets are far from the sun and have a moon,
which makes this a valid attribute implication.

The  implication with $\text{ID}\,  2$ states: All large objects are also far from the sun and have a moon.

The  implication with $\text{ID}\,  3$ states: All objects near to the sun are also small.

The  implication with $\text{ID}\,  4$ states: All  objects far from the sun also have a moon.

The  implication with $\text{ID}\,  5$ states: All  objects without moon are also near the sun and small.

\vfill
\subsection{Attribute exploration}\label{secConExpAttributeExploration}
Implications can also be used for a step-wise computer-guided construction of conceptual knowledge.
The process is
called \emph{attribute exploration} and described in the book by Ganter \& Obiedkov \cite{ganter2016conceptual}.
The developers of \ConExp\ made attribute exploration  available behind the button
{%
\setlength{\fboxsep}{0pt}%
\setlength{\fboxrule}{0.25pt}%
\raisebox{-0.2\height}{\fbox{\includegraphics[scale=0.5]{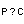}}}%
}
and it can thus be explored by users.
This functionality consists of automatically generated questions which ask if implications are true for all objects.
The user can decide to answer ``no'' and \ConExp\ asks to add an object that represents a counter example to the implication.
This is fairly self-explanatory to follow and is explained in the mentioned documentation and user guide.

\clearpage

\section{Starting \ConExp\ using \Oracle\ \Java\ versions}\label{secInstallation}

\subsubsection*{Problem}
The software \ConExp\  was written more than 18 years ago using a 32\,bit \textsc{Sun} \Java\ version.
Today, \textsc{Sun} \Java\ and compatible newer (\Java\ SE) versions belong to \Oracle,
and they can be downloded from the 
\href{https://www.oracle.com/java/technologies/downloads/archive}{\Oracle\ \Java\ web-archive}.
These \Java\ SE versions are mostly closed-source, and thus open source \Java\ versions may not be completely compatible.
In fact, the standard \Java\ version on Linux is open-source 
and if we start \ConExp\ using this open source \Java, then \ConExp\ will \textbf{NOT  draw concept lattices}.

\subsubsection*{Desire}
Our objective is to have  \ConExp\ draw concept lattices, again.

\subsubsection*{Solution}
A compatible \Oracle\ \Java\ version is needed, some are listed below.
Below, we describe the installation steps,
that 
allowed us to make  \ConExp\ work, \emph{i.e.}, to draw lattices,  on a quite modern 64-bit Linux, which is, Ubuntu 22.04.
The steps do not change  on Ubuntu 20.04,
in fact, they  should also work (with minor adjustments) on other Debian based Linux distributions.

In order to avoid confusion during the installation, we fix a location (in the home directory) where \Java\ and \ConExp\ will be installed or extracted.

\smallskip
\subsection{Location of \ConExp}\label{secHomeConExp}
We use the folder:  
\codebckg{\$HOME/Documents/conexp/}%
,  it has to be created.

\medskip
\subsection{Downloading \ConExp}\label{secInstallConExp}
\providecommand{\README}{\textsc{Readme}}
We  download  
\codebckg{conexp-1.3.zip}
from 
\begin{center}
\url{https://sourceforge.net/projects/conexp/}
\end{center}
and extract it into the sub-folder 
\codebckg{\$HOME/Documents/conexp/conexp-1.3/}.
Among the extracted files there is a \codebckg{readme.txt},
which describes the program and gives instructions.
It is natural, 18 years after its publication,
that following these instructions imposes some difficulties.
Most importantly,  the recommended \textsc{Sun} \Java\ version \codebckg{1.4.1\_01} is not anymore available, not even in the
\href{https://www.oracle.com/java/technologies/downloads/archive}{web-archive}.
Moreover, if we start  \ConExp\ using a default \Java\ version on the system, as described with 
\codebckg{java -jar conexp.jar}%
\footnote{For this we have to be in the folder: \codebckg{\$HOME/Documents/conexp/conexp-1.3/}.},
then the software starts and most functionalities of Section~\ref{secConExp} work, 
but \ConExp\ will \textbf{NOT draw concept lattices}
as shown in   Section~\ref{secConExpBuildLattice}.

\medskip
\GoodNews\ 
The 
\href{https://www.oracle.com/java/technologies/downloads/archive/}{\Oracle\ \Java\ web-archive}
 still contains  compatible  versions, for instance:

\begin{center}
\begin{tabular}[h]{llccl}
		\toprule
		\Java\ version & available for &  32\,bit &  64\,bit & \small verified on Ubuntu 20.04, 22.04 with: \\
		\otoprule
		\href{https://www.oracle.com/java/technologies/java-archive-javase-v14-downloads.html}{Java SE 1.4} 
		& Linux, Windows      &  \green     &         
		& \yellow[\texttt{j2re-1\_4\_2\_19-linux-i586.bin}]
		\\
		\href{https://www.oracle.com/java/technologies/java-archive-javase5-downloads.html}{Java SE 5}   
		& Linux, Windows      &  \checkmark & \green  
		& \green[\texttt{jre-1\_5\_0\_22-linux-amd64.bin}]
		\\
		\href{https://www.oracle.com/java/technologies/javase-java-archive-javase6-downloads.html}{Java SE 6}
		& Linux, Windows      &  \checkmark & \green  
		& \green[\texttt{jre-6u45-linux-x64.bin}]
		\\
		\href{https://www.oracle.com/java/technologies/javase/javase7-archive-downloads.html}{Java SE 7}
		& Linux, MAC, Windows &  \checkmark & \green 
		& \green[\texttt{jre-7u80-linux-x64.tar.gz}]
		\\
		\href{https://www.oracle.com/java/technologies/javase/javase9-archive-downloads.html}{Java SE 9}
		& Linux, MAC, Windows &  \checkmark & \green  
		& \green[\texttt{jre-9.0.4\_linux-x64\_bin.tar.gz}]
		\\
		\bottomrule
\end{tabular}
\end{center}

\medskip
We checked that \ConExp\ draws lattices, on Ubuntu 20.04 and 22.04,  using each of the \Java\ versions shown in the table.
Since these  versions are also available on Windows and MAC, 
we expect  that \ConExp\ works on both.
We describe in Section~\ref{secInstallJava64}
how to install one of the green highlighted 64\,bit versions, which is less complicated on our 64\,bit Linux.
Moreover, we show in Section~\ref{secInstallJava}  how to start \ConExp\ using  the yellow highlighted 32\,bit version, which is more challenging.

\clearpage
\subsection{Get and install  \Java\ SE 5}\label{secInstallJava64}
First, we download the last update of this 
\Java\ SE Runtime Environment, \emph{i.e.}, the version \codebckg{1.5.0\_22}.
A self-extracting archive that installs the second version can be downloaded from the web-archive:
\begin{center}
		\href{https://www.oracle.com/java/technologies/java-archive-javase5-downloads.html}{\texttt{www.oracle.com/java/technologies/java-archive-javase5-downloads.html}}
\end{center}
This web-site contains many \Java\ versions in many sections,
the most recent one is:
\begin{center}
	\includegraphics[width=1\textwidth]{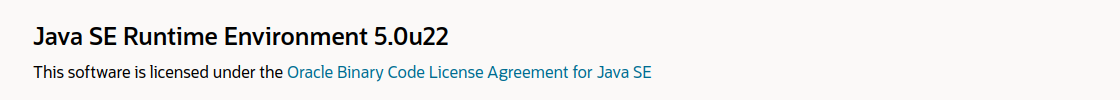}
\end{center}
In this section there are self-extracting files for Windows, Solaris and Linux.
When using a Debian-based 64\,bit Linux (like Ubuntu 20.04 or 22.04) 
the most appropriate self-extracting file is:
\begin{center}
	\includegraphics[width=1\textwidth]{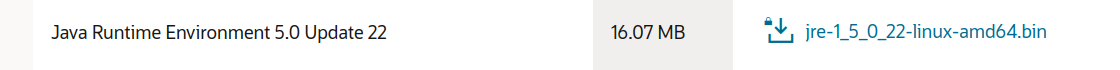}
\end{center}

We  download the file~\codebckg{jre-1\_5\_0\_22-linux-amd64.bin} and save it to \codebckg{\$HOME/Documents/conexp/}.
Next, we have to run it.
For this we need to make it executable.
If the  current working directory of the terminal is
\codebckg{\$HOME/Documents/conexp/},
then the complete process consists of the following 2 instructions:\par
\codebckg{chmod u+x jre-1\_5\_0\_22-linux-amd64.bin}\par
\codebckg{./jre-1\_5\_0\_22-linux-amd64.bin}\\
After agreeing to the licensing terms, \Java\ will be extracted to a sub-folder called \codebckg{./jre1.5.0\_22}.

\bigskip
\subsection{Starting \ConExp}\label{secStartingConExp}
If you completed the steps in Sections~\ref{secHomeConExp}, \ref{secInstallConExp} and  \ref{secInstallJava64}, 
and you put the command 
\codebckg{./jre1.5.0\_22/bin/java -jar conexp-1.3/conexp.jar}
into an executable shell-script into the folder \codebckg{\$HOME/Documents/conexp/} named \codebckg{runConexp.sh},
then 
 running \codebckg{./runConexp.sh} in the terminal
should start  \ConExp.
After performing the steps shown in Sections~\ref{secConExpContext}, 
\ref{secConExpImplications} 
and
\ref{secConExpBuildLattice} 
you will have a window looking as follows:
\begin{center}
\includegraphics[scale=0.3]{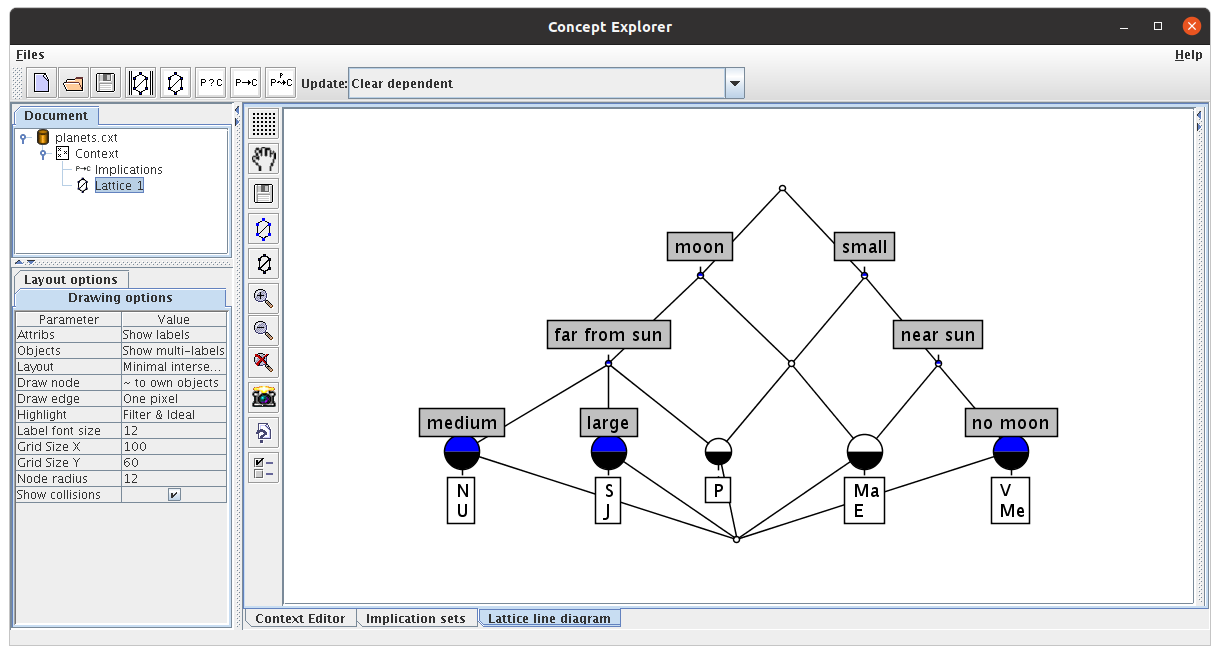}
\end{center}
This window looks a little more modern than our screenshots in Section~\ref{secConExpBuildLattice}, which were taken using the 32\,bit version.

\clearpage
\section{Starting \ConExp\ using a 32  bit \Oracle\ \Java\ version}\label{secInstallation32}

This section is mainly devoted to the installation of a 32\,bit \Java\ version
that allows starting \ConExp\ with the capability of  drawing concept lattices on a 64\,bit Linux.
A few steps are the same:

\begin{itemize}\setlength{\itemsep}{1ex}
\item Just as in Section~\ref{secHomeConExp}, we create the folder \codebckg{\$HOME/Documents/conexp/}.
\item As in Section~\ref{secInstallConExp},
we  download 
\codebckg{conexp-1.3.zip}
from 
\begin{center}
\url{https://sourceforge.net/projects/conexp/}
\end{center}
and extract it into the sub-folder 
\codebckg{\$HOME/Documents/conexp/conexp-1.3/}.

\end{itemize}

\subsection{Get and install  \Java\ SE 1.4}\label{secInstallJava}
We found that \ConExp\ can still draw lattices if we start it using the 
\Java\ SE Runtime Environment of version \codebckg{1.4.1\_07} or \codebckg{1.4.2\_19}.
A self-extracting archive that installs the second version can be downloaded from the web-archive:
\begin{center}
	\href{https://www.oracle.com/java/technologies/java-archive-javase5-downloads.html}{\texttt{www.oracle.com/java/technologies/java-archive-javase-v14-downloads.html}}
\end{center}
This web-site contains many \Java\ versions in many sections,
the most recent one is:
\begin{center}
	\includegraphics[width=\textwidth]{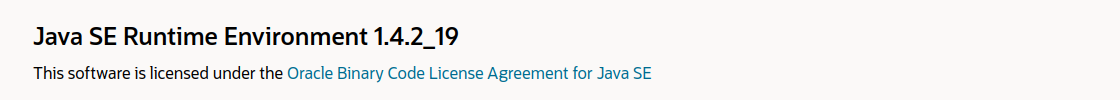}
\end{center}
This section does not contain 64\,bit versions for Linux.
Hence, we use:
\begin{center}
	\includegraphics[width=\textwidth]{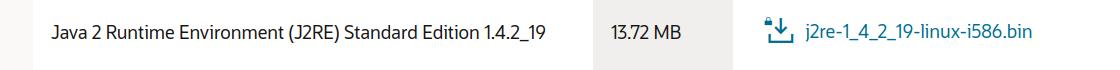}
\end{center}
We  download the file~\codebckg{j2re-1\_4\_2\_19-linux-i586.bin} and save it to \codebckg{\$HOME/Documents/conexp/}.
Next, we have to run it.
However, in our  first attempt  the self-extraction failed, because a 32\,bit ``zip'' library was not installed.
The self-extraction worked after we installed \codebckg{lib32z1}.
If the  current working directory of the terminal is
\codebckg{\$HOME/Documents/conexp/},
then the complete process consists of the following 3 instructions:\par
\codebckg{sudo apt install lib32z1}\par
\codebckg{chmod u+x j2re-1\_4\_2\_19-linux-i586.bin}\par
\codebckg{./j2re-1\_4\_2\_19-linux-i586.bin}\\
After agreeing to the licensing terms, \Java\ will be extracted to a sub-folder called \codebckg{./j2re1.4.2\_19}.

\medskip
\subsection{Installing missing libraries}\label{secInstallLibs}
As of this point, in a terminal with the current working directory \codebckg{\$HOME/Documents/conexp/},
we can run: 
\codebckg{j2re1.4.2\_19/bin/java -jar conexp-1.3/conexp.jar}.

The first run fails with the \Java\ Exception:
\begin{verbatim}
    Exception in thread "main" java.lang.UnsatisfiedLinkError:
    $HOME/Documents/conexp/j2re1.4.2_19/lib/i386/libawt.so: libXt.so.6: 
    cannot open shared object file: No such file or directory
\end{verbatim}
This exception will be avoided after  installing the library \codebckg{libxt6:i386}.
On the second and third run similar exceptions appear and each time they can be ``fixed'' by installing a library.
We had to install 3 libraries (including dependencies)  and a font package using the following 4 commands:\par
\codebckg{sudo apt install libxt6:i386}\par
\codebckg{sudo apt install libxext6:i386}\par
\codebckg{sudo apt install libxtst6:i386}\par
\codebckg{sudo apt install xfonts-100dpi}  \hfill(These fonts remove a font-warning, after a reboot.)

At this point,
in a terminal with the current working directory \codebckg{\$HOME/Documents/conexp/},
running 
\codebckg{j2re1.4.2\_19/bin/java -jar conexp-1.3/conexp.jar}
will produce  an almost zero-width,  \textbf{barely visible}  ``Concept Explorer'' window,
which looks as follows:
\begin{center}
	\includegraphics{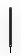}	
\end{center}
Once found, the window can be pulled to the desired size.
This  \Java\ bug
 disappears as of version Java SE~5.
\GoodNews\ 
There is also a  work-around for this bug.
The thin ``Concept Explorer'' window can be found, resized and put in front automatically
from the terminal
using the command \codebckg{wmctrl},
which can to be installed using the command: 
\codebckg{sudo apt install wmctrl}.

\medskip
\subsection{Starting \ConExp\ with a 32 bit \Java}\label{secStartingConExp32}
Once the libraries, the font and  \codebckg{wmctrl} are installed,
 \ConExp\ can be started running  the following shell-script that
 collects 4 commands and should be saved in the folder  \codebckg{\$HOME/Documents/conexp/}.

\begin{flushleft}
\begin{Verbatim}
#!/bin/sh
./j2re1.4.2_19/bin/java -jar conexp-1.3/conexp.jar  &
sleep 2
wmctrl -r 'Concept Explorer' -e 0,0,0,900,600
wmctrl -a 'Concept Explorer'
\end{Verbatim}
\end{flushleft}
 Running the script \codebckg{./runConexp.sh} in the terminal
should start  \ConExp, position the ``Concept Explorer'' window into the top left corner of the screen, resize it  to $900\times 600$ pixels and activate it.
The window should look like in Section~\ref{secConExpContext}.

If the window is still thin, then you may need to increase the sleep timer.
The name of the \ConExp\ window should be  ``Concept Explorer'',
you can check it with \codebckg{wmctrl -l}.

We hope that the instructions above simplify the installation and the use of the Concept Explorer
to the reader.

\vfill


\end{document}